\documentclass[11pt]{article}
\usepackage{amscd}
\usepackage{amsfonts}
\usepackage{amsmath}
\usepackage{amssymb}
\usepackage{amsthm}
\usepackage{bbm}
\usepackage{CJK}
\usepackage{fancyhdr}
\usepackage{graphicx}
\usepackage{hyperref}
\usepackage{indentfirst}
\usepackage{latexsym}
\usepackage{mathrsfs}
\usepackage{xypic}

\usepackage[top=1in,bottom=1in,left=1.25in,right=1.25in]{geometry}
\def\<{\langle}
\def\>{\rangle}
\def\a{\alpha}
\def\b{\beta}

\def\c{\cdot}

\def\D{\Delta}

\def\g{\gamma}

\def\o{\otimes}

\def\v{\varepsilon}

\begin{document}
\renewcommand{\baselinestretch}{1.2}
\renewcommand{\arraystretch}{1.0}
\title{\textbf{On the structure theorem and the Maschke type theorem of Doi Hom-Hopf modules  }}
\author{{Shuang-jian Guo$^{1}$\thanks{Corresponding author, email: shuangjguo@gmail.com } \quad   Xiao-hui Zhang$^{2}$}\\
{\small 1. School of Mathematics and Statistics, Guizhou University of Finance and Economics}\\
{\small Guiyang,  550025, P. R. China}\\
{\small 2. Department of Mathematics, Southeast  University}\\
 {\small Nanjing, 210096,  P. R. China}}
 \date{}
 \maketitle
\noindent\textbf {\bf Abstract}   We give necessary and sufficient conditions for the functor that forgets the $(C, \gamma)$-coaction to be separable. This leads to a generalized notion of integrals. Finally,  the applications of our results are considered.

\noindent\textbf {Keywords } Monoidal Hom-Hopf algebra; Doi Hom-Hopf module; separable functors;    normalized $(A,\b)$-integral;   Maschke type theorem.

\noindent\textbf {MR(2010)Subject Classification} 16T05

\section{Introduction}
The present paper investigates variations on the theme of
Hom-algebras, a topic which has recently received much attention
from various researchers. The study of Hom-associative algebras
originates with work by Hartwig, Larsson and Silvestrov in the Lie
case \cite{HLS}, where a notion of Hom-Lie algebra was introduced in the
context of studying deformations of Witt and Virasoro algebras.
Later, it was extended to the associative case by Makhlouf and
Silverstrov in \cite{MS08}. Now the associativity is replaced by
Hom-associativity $\alpha(a)(bc)=(ab)\alpha(c)$. Hom-coassociativity
for a Hom-coalgebra can be considered in a similar way, see \cite{MS08}.
Caenepeel etc. \cite{CG11} studied Hom-structures from the point of view of
monoidal categories. This leads to the natural definition of monoidal Hom-algebras,
Hom-coalgebras, etc. They  constructed a symmetric monoidal category, and then
introduced monoidal Hom-algebras, Hom-coalgebras, etc. as algebras, coalgebras,
etc. in this monoidal category.

Let $H$ be a Hopf algebra with bijective antipode over a commutative ring $k$. As
is well-known, a Hopf module is a $k$-module that is at once an $H$-module and an
$H$-comodule, with a certain compatibility relation; see \cite{S69} for details. Doi \cite{D92}
generalized this concept in the following way: if $A$ is an $H$-comodule algebra, and
$C$ is an $H$-module coalgebra, then he introduced a so-called unified Hopf module:
this is a $k$-module that is at once an $A$-module and a $C$-comodule, satisfying a
compatibility relation that is an immediate generalization of the one that may
be found in Sweedler's book \cite{S69}. One of the nice features here is that Doi's
Hopf modules (we will call them Doi-Hopf modules) really unify a lot of module
structures that have been studied by several authors; let us mention Sweedler's
Hopf modules \cite{S69}, Takeuchi's relative Hopf modules \cite{T}, graded modules, and
modules graded by a $G$-set. In \cite{D92}, induction functors between categories of Doi-Hopf modules and their adjoints are studied, and it turns out that many pairs
of adjoint functors studied in the literature (the forgetful functor and its adjoint,
extension and restriction of scalars,...) are special cases.  Caenepeel etc. \cite{CGBZ99} proved a Maschke type theorem for the category of Doi-Hopf modules. In fact, they gave necessary and sufficient conditions for the functor that forgets the $H$-coaction to be separable. This leads to a generalized notion of integrals of Doi \cite{D85}.

The following questions arise naturally:

1. How do we introduce the notion of a Doi Hom-Hopf module?

2. The evident functor $F$ from the category
of Doi Hom-Hopf modules to the category of modules over the Hom algebra $A$, does it possess a right adjoint? What is a sufficient and necessary condition for its
separability?

3.~ How do we give a Maschke type theorem for Doi Hom-Hopf modules? The aim of this article is to give a positive answer to these
questions.

In this paper we study the generalization of the previous results to the monoidal Hom-Hopf algebras.
In Sec.3, we introduce the notion of a Doi Hom-Hopf module and prove that the functor $F $ from  the category of Doi Hom-Hopf modules to the category of right $(A, \b)$-Hom-modules has a right
 adjoint(see Proposition 3.4).
In Sec.4, we obtain the main result of this paper.
 we give necessary and sufficient conditions for the functor that forgets the $(C, \g)$-coaction to be separable (see Theorem 4.2).  The applications of our results are considered in Sec.5.

\section{Preliminaries}
Throughout this paper we work over a commutative ring $k$, we recall from \cite{CG11}
some information about Hom-structures which are needed in what follows.

Let $\mathcal{C}$ be a category. We introduce a new category
$\widetilde{\mathscr{H}}(\mathcal{C})$ as follows: objects are
couples $(M, \mu)$, with $M \in \mathcal{C}$ and $\mu \in
Aut_{\mathcal{C}}(M)$. A morphism $f: (M, \mu)\rightarrow (N, \nu )$
is a morphism $f : M\rightarrow N$ in $\mathcal{C}$ such that $\nu
\circ f= f \circ \mu$.

Let $\mathscr{M}_k$ denotes the category of $k$-modules.
~$\mathscr{H}(\mathscr{M}_k)$ will be called the Hom-category
associated to $\mathscr{M}_k$. If $(M,\mu) \in \mathscr{M}_k$, then
$\mu: M\rightarrow M$ is obviously a morphism in
~$\mathscr{H}(\mathscr{M}_k)$. It is easy to show that
~$\widetilde{\mathscr{H}}(\mathscr{M}_k)$ =
(~$\mathscr{H}(\mathscr{M}_k),~\otimes,~(I, I),~\widetilde{a},
~\widetilde{l},~\widetilde{r}))$ is a monoidal category by
Proposition 1.1 in \cite{CG11}: the tensor product of $(M,\mu)$ and $(N,
\nu)$ in ~$\widetilde{\mathscr{H}}(\mathscr{M}_k)$ is given by the
formula $(M, \mu)\otimes (N, \nu) = (M\otimes N, \mu \otimes \nu)$.

Assume that $(M, \mu), (N, \nu), (P,\pi)\in
\widetilde{\mathscr{H}}(\mathscr{M}_k)$. The associativity and unit
constraints are given by the formulas
\begin{eqnarray*}
\widetilde{a}_{M,N,P}((m\o n)\o p)=\mu(m)\o (n\o \pi^{-1}(p)),\\
\widetilde{l}_{M}(x\o m)=\widetilde{r}_{M}(m\o x)=x\mu(m).
\end{eqnarray*}
An algebra in $\widetilde{\mathscr{H}}(\mathscr{M}_k)$ will be called a monoidal Hom-algebra.\\

\noindent{\bf Definition 2.1.} A monoidal Hom-algebra  is an
object is a triple $(A, \alpha)\in
\widetilde{\mathscr{H}}(\mathscr{M}_k)$ together with a $k$-linear
map $m_A: A\o A\rightarrow A$ and an element $1_A\in A$ such that
\begin{eqnarray*}
\alpha(ab)=\alpha(a)\alpha(b); ~\a(1_A)=1_A,\\
\alpha(a)(bc)=(ab)\alpha(c); ~a1_A=1_Aa=\a(a),
\end{eqnarray*}
for all $a,b,c\in A$. Here we use the notation $m_A(a\o b)=ab$.

\noindent{\bf Definition 2.2.} A monoidal Hom-coalgebra  is
an object $(C,\g)\in \widetilde{\mathscr{H}}(\mathscr{M}_k)$
together with $k$-linear maps $\Delta:C\rightarrow C\otimes C,~~~
\D(c)=c_{(1)}\o c_{(2)}$ (summation implicitly understood) and
$\gamma:C\rightarrow C$ such that
\begin{eqnarray*}
\D(\g(c))=\g(c_{(1)})\o \g(c_{(2)}); ~~~\varepsilon(\g(c))=\varepsilon(c),
\end{eqnarray*}
and
\begin{eqnarray*}
\g^{-1}(c_{(1)})\otimes c_{(2)(1)}\otimes c_{(2)(2)}=c_{(1)(1)}\otimes c_{(1)(2)}\otimes\g^{-1}(c_{(2)}),~\varepsilon(c_{(1)})c_{(2)}=\varepsilon(c_{(2)})c_{(1)}=\g^{-1}(c)
\end{eqnarray*}
for all $c\in C$.\\

\noindent{\bf Definition 2.3.} A monoidal Hom-bialgebra
$H=(H,\a,m,\eta, \Delta,\varepsilon)$  is a bialgebra in
the symmetric monoidal category
$\widetilde{\mathscr{H}}(\mathscr{M}_k)$. This means that $(H, \a,
m,\eta)$ is a Hom-algebra, $(H,\Delta,\alpha)$ is a Hom-coalgebra
and that $\D$ and $\v$ are morphisms of Hom-algebras that is,
\begin{eqnarray*}
\Delta(ab)=a_{(1)}b_{(1)}\otimes a_{(2)}b_{(2)}; ~~~\Delta(1_H)=1_H\otimes 1_H,\\
\varepsilon(ab)=\varepsilon(a)\varepsilon(b),~\varepsilon(1_H)=1_H.
\end{eqnarray*}

\noindent{\bf Definition 2.4.} A monoidal Hom-Hopf algebra is
a monoidal Hom-bialgebra $(H, \alpha)$ together with a linear map
$S:H\rightarrow H$ in $\widetilde{\mathscr{H}}(\mathscr{M}_k)$ such
that
$$S\ast I=I\ast S=\eta\varepsilon,~S\alpha=\alpha S.$$

\noindent{\bf Definition 2.5.} Let $(A,\alpha)$ be a monoidal
Hom-algebra.  A right $(A,\alpha)$- Hom-module is an object
$(M,\mu)\in \widetilde{\mathscr{H}}(\mathscr{M}_k) $ consists of a
$k$-module and a linear map $\mu:M\rightarrow M$ together with a
morphism $\psi:M\otimes A\rightarrow M, \psi(m\c a)=m\c a$, in
$\widetilde{\mathscr{H}}(\mathscr{M}_k) $  such that
\begin{eqnarray*}
(m\c a)\c\alpha(b)=\mu(m)\c (ab);~~~~ m\c 1_A=\mu(m),
\end{eqnarray*}
for all $a\in A$ and $m\in M$. The fact that $\psi\in \widetilde{\mathscr{H}}(\mathscr{M}_k)$ means that
\begin{eqnarray*}
\mu(m\c a)=\mu(m)\c\alpha(a).
\end{eqnarray*}
A morphism $f: (M, \mu)\rightarrow (N, \nu)$ in $\widetilde{\mathscr{H}}(\mathscr{M}_k)$ is called right $A$-linear if it preserves the $A$-action, that is, $f(m\c a)=f(m)\c a$. $\widetilde{\mathscr{H}}(\mathscr{M}_k)_A$ will denote the category of right $(A, \a)$-Hom-modules and $A$-linear morphisms.\\

\noindent{\bf Definition 2.6.} Let $(C,\g)$ be a monoidal
Hom-coalgebra. A right $(C,\g)$-Hom-comodule  is an object
$(M,\mu)\in \widetilde{\mathscr{H}}(\mathscr{M}_k)$  together with a
$k$-linear map $\rho_M: M\rightarrow M\otimes C$ notation
$\rho_M(m)=m_{[0]}\o m_{[1]}$ in
$\widetilde{\mathscr{H}}(\mathscr{M}_k)$ such that
\begin{eqnarray*}
m_{[0][0]}\otimes (m_{[0][1]}\otimes
\g^{-1}(m_{[1]}))=\mu^{-1}(m_{[0]})\otimes \D_C(m_{[1]});
~m_{[0]}\varepsilon(m_{[1]})=\mu^{-1}(m),
\end{eqnarray*}
for all $m\in M$.  The fact that $\rho_M\in \widetilde{\mathscr{H}}(\mathscr{M}_k)$ means that
\begin{eqnarray*}
\rho_M(\mu(m))=\mu(m_{[0]})\otimes \g(m_{[1]}).
\end{eqnarray*}
Morphisms of right $(C, \g)$-Hom-comodule are defined in the obvious way. The
category of right $(C, \g)$-Hom-comodules will be denoted by $\widetilde{\mathscr{H}}(\mathscr{M}_k)^C$.
\\

\noindent{\bf Theorem 2.7.} {\bf (Rafael Theorem)}
Let $L:\mathcal{C}\rightarrow \mathcal{D}$ be the left adjoint functor of $R:\mathcal{D}\rightarrow \mathcal{C}$. Then $L$ is a separable functor if and only if the unit $\eta$ of the adjunction $(L,R)$ has a natural retraction, i.e., there is a natural transformation $\nu: RL\rightarrow id_{\mathcal{C}}$, such that
$\nu \circ \eta= id$.

\section{Adjoint functor}

\noindent{\bf Definition 3.1.} Let $(H,\alpha)$  be a monoidal Hom-Hopf algebra.  A monoidal Hom-algebra  $(A,\beta)$ is called  a right
$(H,\alpha)$-Hom-comodule algebra if $(A,\beta)$ is a right $(H,\alpha)$ Hom-comodule with coaction $\rho_A: A\rightarrow A\o H,
 ~~\rho_A(a)= a_{[0]}\o a_{[1]}$
 such that the following conditions satisfy,
 \begin{eqnarray*}
 &&\rho_{A} ( ab)= a_{[0]} b_{[0]} \otimes a_{[1]}b_{[1]},\\
 &&\rho_{A} ( 1_A)=1_A\o 1_H.
 \end{eqnarray*}
for all $a,b \in A$.\\

\noindent{\bf Definition 3.2.} Let $(H,\alpha)$ be a monoidal Hom-Hopf algebra. A monoidal Hom-coalgebra $(C, \g)$ is called a right $(H,\alpha)$ Hom-module coalgebra, if $(C,\g)$ is a right $(H,\alpha)$ Hom-module with action $\phi: C \otimes H\rightarrow C$,
$\phi (c \otimes h) = c \cdot h$ such that the following conditions satisfy:
\begin{eqnarray*}
&&\Delta(c \cdot h)= c _{(1)}\cdot h _{(1)} \otimes c _{(2)}\cdot h _{(2)},\\
&&\varepsilon(c \cdot h) = \varepsilon(c)\varepsilon(h),
\end{eqnarray*}
for all  $c \in C$ and $g, h \in H$.

A  Doi Hom-Hopf datum is a triple ($H, A, C$), where $H$ is a monoidal Hom-Hopf algebra, $A$
a  right $(H, \a)$-Hom comodule algebra and $(C, \g)$ a  right  $(H, \a)$-Hom module coalgebra. \\

\noindent{\bf Definition 3.3.}  Given a  Doi Hom-Hopf datum ($H, A, C$).  A Doi Hom-Hopf module $(M,\mu)$
is a right $(A,\b)$-Hom-module which is also a right
$(C,\g)$-Hom-comodule with the  coaction structure $ \rho _{M} : M
\rightarrow M \otimes  C $ defined by $\rho _{M}(m)=m_{[0]}\o
m_{[1]}$ such that the following compatible condition holds: for all
$m\in M $ and $a\in A$,
\begin{eqnarray*}
&&\rho _{M}( m\c a) = m _{[0]} \c a _{[0]} \otimes  m _{[1]}a _{[1]}.
\end{eqnarray*}

A morphism between two right Doi Hom-Hopf modules is a $k$-linear map
which is a morphism in the categories $\widetilde{\mathscr{H}}(\mathscr{M}_k)_A$ and $\widetilde{\mathscr{C}}(\mathscr{M}_k)^C$ at the same time.  $ \widetilde{\mathscr{H}}(\mathscr{M}_k)(H)^{C}_{A}$ will denote  the category of right Doi Hom-Hopf modules and morphisms
between them. \\

\noindent{\bf Proposition 3.4.}
 The forgetful functor
$F:  \widetilde{\mathscr{H}}(\mathscr{M}_k)(H)^{C}_{A} \rightarrow \widetilde{\mathscr{H}}(\mathscr{M}_k)_{A}$ has a right
 adjoint $G:  \widetilde{\mathscr{H}}(\mathscr{M}_k)_{A}\rightarrow \widetilde{\mathscr{H}}(\mathscr{M}_k)(H)^{C}_{A}$.
 $G$ is defined by
$$
 G(M) = M \otimes C,
 $$
 with structure maps
\begin{eqnarray*}
( m  \otimes c  ) \cdot a =   m \c a _{[0]}\otimes c \c a _{[1]},
\end{eqnarray*}
\begin{eqnarray*}
\rho _{G(M)} ( m  \otimes c ) =  (\mu^{-1}( m ) \otimes c _{(1)}) \otimes \g(c _{(2)}),
\end{eqnarray*}
for all $a \in A $ and $m\in M, c\in C$.

{\it Proof.}
Let us first show that $G (M)$ is an object of $ \widetilde{\mathscr{H}}(\mathscr{M}_k)(H)^{C}_{A}$ . It is routine to check that $G (M)$ is a  right $(C,\g)$-Hom-comodule and a right $(A,\b)$-Hom-module. Now we only check the compatibility condition, for all $a \in A $ and $c\in C$.
 Indeed,
\begin{eqnarray*}
 \rho_{G(M)}  ((m  \otimes c)\cdot a) & =&  \rho_{G(M)}  ( m \c a _{[0]} \otimes c \c a _{[1]})\\
 & = & \mu^{-1}(m)\c \b^{-1}(a_{[0]})\o c_{(1)} \c a _{[1](1)}\o \gamma(c_{(2)}\c a _{[1](2)})\\
 & = & \mu^{-1}(m)\c a_{[0][0]}\o c_{(1)}\c a _{[0][1]}\o \gamma(c_{(2)})\c a _{[1]}\\
 & = & (m  \otimes c)_{[0]}\c a_{[0]}\o (m  \otimes c)_{[1]}\c a _{[1]}\\
 &=& \rho_{G(M)} (m\otimes c )\cdot a.
  \end{eqnarray*}
This is exactly what we have to show.

For an $A$-linear map $\varphi : (M, \mu) \rightarrow (N, \nu)$, we put
 $$
 G (\varphi)  = \varphi \otimes id _{C}:
M \otimes C\rightarrow N\otimes C .
 $$
Standard computations show that $G (\varphi)$ is  morphisms of right $(A,\b)$-Hom-modules and right $(C,\g)$-Hom-comodules. Let us describe the
 unit $\eta$ and the counit $\delta$ of the adjunction. The unit is described by the coaction:
 for $(M,\mu) \in   \widetilde{\mathscr{H}}(\mathscr{M}_k)(H)^{C}_{A} $,  we define
 $\eta _{M} :
  M \rightarrow M \otimes C$  as follows:
 for all $m \in M $,
 $$
 \eta _{M}(m) =  m _{[0]}
 \otimes m _{[1]}.
 $$
 We can check that $\eta_{M} \in  \widetilde{\mathscr{H}}(\mathscr{M}_k)(H)^{C}_{A} $ . In fact, for any $m\in M$, we have
 \begin{eqnarray*}
  \eta _{M}(m\c a)&=& (m\c a)_{[0]}\o (m\c a)_{[1]}\\
  &=& m_{[0]}\c a_{[0]}\o m_{[1]}a_{[1]}\\
  &=& (m_{[0]}\o m_{[1]})\c a=\eta _{M}(m)\c a,
 \end{eqnarray*}
and
\begin{eqnarray*}
\rho_{M\o C}\circ   \eta _{M}(m)&=& \rho_{M\o C}(m _{[0]}
 \otimes m _{[1]})\\
 &=&(\mu^{-1}( m_{[0]} ) \otimes m _{[1](1)}) \otimes \g(m _{[1](2)})\\
 &=&( m_{[0][0]}  \otimes m _{[0][1]}) \otimes m _{[1]}\\
&=& ( \eta _{M}\o id_C)(m_{[0]}\o m_{[1]})\\
 &=&( \eta _{M}\o id_C)\circ \rho_{M}(m).
\end{eqnarray*}
  For any $(N, \nu) \in \widetilde{\mathscr{H}}(\mathscr{M}_k)_{A}$, we define
 $\delta _{N}: N\otimes C\rightarrow N$, for all $n \in N$ and $c \in C $,
 $$
 \delta _{N}(n  \otimes c  )
 = \varepsilon(c ) \nu(n) ,
 $$
 we can check that $\delta _{N}$ is $(A,\b)$-linear.  In fact, for any $n\in N$, we have
\begin{eqnarray*}
 \delta _{N}((n  \otimes c )\c a )&=&  \delta _{N}(n\c a_{[0]}  \otimes c\c a_{[1]}  )\\
 &=&\varepsilon(c\c a_{[1]})\nu(n\c a_{[0]})\\
 &=&\varepsilon(c)\nu(n)\c a=  \delta _{N}(n  \otimes c  )\c a.
\end{eqnarray*}
This is what we need to show. We can check that $\eta$ and $\delta$ defined above are all natural transformations and
 satisfied
 $$
 G (\delta _{N})\circ \eta _{G (N)} = I _{G (N)},
 $$
 $$
 \delta _{F (M)}\circ F (\eta _{M}) = I _{F (M)},
 $$
 for all $(M,\mu) \in \widetilde{\mathscr{H}}(\mathscr{M}_k)(H)^{C}_{A} $ and $(N,\nu) \in  \widetilde{\mathscr{H}}(\mathscr{M}_k)_{A}$ .

 \section{\textbf{Separable functors for the category of Doi Hom-Hopf modules }}
\def\theequation{4. \arabic{equation}}
\setcounter{equation} {0} \hskip\parindent

In this section, we shall  give necessary and sufficient conditions for the functor $F$
 which forget the $(C, \g)$-coaction to be separable.\\

\noindent{\bf Definition 4.1.} Let ($H, A, C$) be a  Doi Hom-Hopf datum.  A $k$-linear map
  $$
\theta : (C, \g)  \otimes (C, \g)  \rightarrow  (A, \b)
  $$ such that $\theta\circ (\g\o \g)=\b\circ \theta$  is called a  normalized $(A, \b)$-integral, if $\theta$ satisfies the following conditions:

  (1) For all $c,d \in C $,
 \begin{eqnarray}
 && \theta (\g^{-1}(d) \otimes c _{(1)} )\otimes \g(c _{(2)}) \nonumber\\
 &=& \beta( \theta  (d _{(2)} \otimes \g^{-1}(c))
   _{[0]})\otimes d _{(1)} \c\theta (d _{(2)} \otimes \g^{-1}(c))_{[1]}.
 \end{eqnarray}

(2) For all  $c \in C $,
 \begin{equation}
 \theta (c _{(1)} \otimes c _{(2)})
 = 1 _{A } \varepsilon (c).
 \end{equation}

(3) For all $a \in A, c,d \in C$,
 \begin{equation}
  \b^2 (a _{[0][0]})
   \theta (\g^{-1}(d) a _{[0][1]} \otimes \g^{-1}(c)\g^{-1}(a _{[1]}))
   = \theta ( d \otimes c ) a.
 \end{equation}

\noindent{\bf Theorem 4.2.}
For any  Doi Hom-Hopf datum ($H, A, C$), the following assertions are equivalent,

(1) The left adjoint functor $F$ in Proposition 3.4 is separable,

(2) There exists a normalized $(A, \b)$-integral
$\theta : (C, \g)  \otimes (C, \g)  \rightarrow  (A, \b) $.

{\it Proof.} $(2) \Longrightarrow(1)$. For any Doi Hom-Hopf module $M$, we define
\begin{eqnarray*}
 \nu_{M} : M\otimes C &\rightarrow& M ,\\
m\otimes c  &\mapsto& \mu(m _{[0]}) \theta  (m _{[1]} \otimes \g^{-1}(c)),
\end{eqnarray*}
for all  $m \in M $ and $c \in C $.
Now, we shall check that $\nu _{M} \in \widetilde{\mathscr{H}}(\mathscr{M}_k)(H)^{C}_{A}$. In fact, for all $m \in M$, $c \in C$ and $a \in A $, it is easy to get that
$$
\nu_M(\mu(m) \o \g(c)) = \mu(\nu_M(m \o c)).
$$
We also have
\begin{eqnarray*}
\nu _{M} ((m  \otimes c )\cdot a)&=&\nu_{M}( m \c a _{[0]} \otimes c \c a _{[1]})\\
&=& (\mu(m_{[0]})\cdot\b (a _{[0][0]})) \theta
    (m_{[1]} a _{[0][1]}\otimes \g^{-1}(c)  \c\a^{-1}(a _{[1]}))\\
&=& \mu^2(m_{[0]})\cdot(\b(a _{[0][0]})\b^{-1}(\theta
    (m_{[1]} a _{[0][1]}\otimes \g^{-1}(c) \c \a^{-1}(a _{[1]})))\\
&=& \mu^2(m_{[0]})\cdot(\b(a _{[0][0]})\theta(\g^{-1}(m_{[1]})\g^{-1}(a _{[0][1]}) \otimes \g^{-2}(c) \c \a^{-2}(a _{[1]})))\\
&\stackrel{(4.3)}{=}& \mu^2(m_{[0]})\cdot(\theta(m_{[1]} \otimes \g^{-1}(c)) \c\b^{-1}(a))\\
&=& (\mu(m_{[0]})\cdot \theta(m_{[1]} \otimes \g^{-1}(c)))\cdot a\\
&=& (\nu _{M}(m  \otimes c ))\cdot a.
\end{eqnarray*}
Hence  it is a morphism of $(A, \b)$-Hom-modules. Next, we shall check that $\nu _{M}$ is a morphism of comodules
over $(C, \g)$. It is sufficient to check that
 $$
 \rho _{M}  \circ \nu_{M}  = (\nu_{M}  \otimes id _{C})
 \circ \rho _{M}
 $$
 holds. For all $m \in M $ and $c \in C$,  we have
\begin{eqnarray*}
&&\rho _{M}\circ \nu_{M}(m  \otimes c )\\
&=& \rho _{M}  (\mu(m_{[0]} )
\theta (m_{[1]} \otimes \g^{-1}(c )  ))\\
&=& (\mu(m_{[0]})\c
\theta(m_{[1]} \otimes \g^{-1}(c))_{[0]}
 \otimes (\mu(m_{[0]})\c \theta (m_{[1]} \otimes \g^{-1}(c) )) _{[1]} \\
&=& \mu(m_{[0][0]})\c \theta (m_{[1]}    \otimes \g^{-1}(c))_{[0]}   \otimes
   \g( m_{[0][1]} )\c
\theta (m_{(1)}\otimes \g^{-1}(c) ) _{[1]} \\
&=& m_{[0]}\c
\theta (\g(m_{[1](2)})\otimes \g^{-1}(c))_{[0]}
 \otimes\g(m_{[1](1)}) \c\theta (\g(m_{[1](2)})\otimes \g^{-1}(c)  ) _{[1]}\\
&\stackrel{(4.1)}{=}& m_{[0]}\c
\b^{-1}(\theta (m_{[1]} \otimes c_{(1)}))
 \otimes \g(c  _{(2)}) \\
& = & m_{[0]}\c
\theta (\g^{-1}(m_{[1]}) \otimes \g^{-1}(c_{(1)}))
 \otimes \g(c_{(2)}) \\
& = &(\nu_{M}  \otimes id _{C})
 \circ \rho _{M} (m \otimes c).
\end{eqnarray*}
For all $m \in M $, since
\begin{eqnarray*}
 \nu _{M}  \circ \eta _{M}  (m) &=& \nu _{M} (m _{[0]}
 \otimes m _{[1]})\\
&= & \mu(m _{[0][0]})\c \theta  (m _{[0][1]}
 \otimes  \g^{-1}(m _{[1]}))\\
&= &  m _{[0]}\c \theta (m _{[1](1)}\otimes m_{[1](2)})\stackrel{(4.2)}{=}
 m.
\end{eqnarray*}
So the left adjoint $F$ in Proposition 3.4 is separable follows by Rafael theorem (see \cite{R90} for detail).

$(1) \Longrightarrow (2)$. We consider the following  Doi Hom-Hopf module $A\o C $,
and the $(A, \b)$-actions and $(C, \g)$-coaction are defined as follows:
$$
\left\{
  \begin{array}{ll}
(a \otimes c ) \cdot b = a b _{[0]}\otimes c\c b _{[1]};\\
  \rho_{A\o C}  ( a \otimes c  )
  =(\b^{-1}(a) \otimes c _{(1)})\otimes \g(c _{(2)}),
\end{array}
\right.
$$
for any $a, b\in A$ and $c\in C$.

Evaluating at this  object, the retraction $\nu$ of the unit of the adjunction in Proposition 3.4 yields a morphism
$$
\nu _{A\o C  }:  (A  \otimes C) \otimes C  \rightarrow A \otimes C
$$
such that, for all  $a \in A , c \in C$,

$$
\nu _{A \o C}((a \otimes c _{(1)})
\otimes c _{(2)})= \b(a) \otimes c.$$

 It can be used to construct $\theta$ as follows:
$$
\theta :C \otimes C\rightarrow A,
$$
$$
\theta (c \otimes d) = r_A (id _{A }\otimes \varepsilon_C) \nu _{A \o C}
((1 _{A }
\otimes c) \otimes d ),
$$
where $r$ means the right unit constraint. For all $c \in C $, since
\begin{eqnarray*}
&&\theta (c _{(1)}\otimes c _{(2)})\\
&=&r_A(id _{A }\otimes \varepsilon_C) \nu _{A \o C}((1 _{A}
\otimes c _{(1)})
\otimes c _{(2)})\\
&=&r_A(id _{A }\otimes \varepsilon_C) (1 _{A}\otimes c )
= 1 _{A} \varepsilon_C(c).
\end{eqnarray*}
Hence condition Eq.(4.2) follows. It can be seen to obey Eq.(4.3) by naturality and the
$(A, \b)$-Hom module map property of $\nu$.

 The verification of Eq.(4.1) is more
involved. For any right $(C, \g)$-comodule $M$, we consider the  Doi Hom-Hopf module
$M \otimes  A$,
the $(A, \b)$-action and $(C,\g)$-coaction are defined as follows: for all $m \in M$ and $a, b \in A $,
$$
\left\{
  \begin{array}{ll}
    (m \otimes a) \cdot b = \mu^{-1}(m) \otimes a\b(b), \\
  \rho_{M\o A} ( m \otimes a)= (m _{[0]} \otimes a _{[0]})\otimes m _{[1]} \c a _{[1]}.
  \end{array}
\right.
$$
In particular, there is a  Doi Hom-Hopf module $C \otimes  A$ and the map
$$
\xi : C \otimes A  \rightarrow A  \otimes C\
$$
$$
\xi(c \otimes a) = \b(a _{[0]}) \otimes \g^{-1}(c) \c a _{[1]}.
$$
Since $\xi$ is both right $(A, \b)$-linear and right $(C, \g)$-colinear, thus we have
\begin{eqnarray}
\xi(\nu_{C \o A}((c \o a) \o d)) &=& \nu_{A \o C}((\xi \o id_C)((c \o a) \o d)) \nonumber\\
&=& \nu_{A \o C}((\b(a_{[0]}) \o \g^{-1}(c)\c a_{[1]})\o d).
\end{eqnarray}
It is not hard to check that $GF(C \o A)= (C\o A) \o C \in {}^{C}\widetilde{\mathscr{H}}(\mathscr{M}_k)^{C}_{A}$, and its left $(C, \g)$-Hom comodule structure
is given by
$$
(c \o a) \o d \mapsto \a(c_{(1)}) \o ((c_{(2)} \o \b^{-1}(a)) \o \g^{-1}(d)).
$$
Also $C \o A\in {}^{C}\widetilde{\mathscr{H}}(\mathscr{M}_k)^{C}_{A}$, and the left $(C,\g)$-coaction of $C \o A$ is given by
$$
c \o a\mapsto \g(c_{(1)}) \o (c_{(2)} \o \b^{-1}(a)).
$$
We also get $\nu_{C\o A}:(C\o A) \o C\rightarrow C \o A$ is a morphism in ${}^{C}\widetilde{\mathscr{H}}(\mathscr{M}_k)^{C}_{A}$, which means
\begin{eqnarray}
&\nu_{C\o A}((c \o a) \o d)_{[-1]} \o \nu_{C\o A}((c \o a) \o d)_{[0]} \nonumber\\
&~~~~~~~~= \g(c_{(1)}) \o \nu_{C\o A}((c_{(2)} \o \b^{-1}(a)) \o \g^{-1}(d)).
\end{eqnarray}
Thus we conclude that $\nu_{C\o A}$ is left and right $(C,\g)$-colinear. Take $c, d\in C$, and put
$$\nu_{A\o C}((1_A \o c) \o d) =\sum_{i} a_i\o q_i \in A \o C,$$
$$\nu_{C\o A}((h\o 1_A)\o g)=\sum_{i} p_i\o b_i \in C \o A,$$
we obtain
\begin{eqnarray*}
&&\b(\theta(c_{(2)} \o \g^{-1}(d))_{[0]}) \o c_{(1)}\theta(c_{(2)} \o \g^{-1}(d))_{[1]}\\
&=& \b(r_A(id_A \o \varepsilon_C)\nu_{A \o C}((1_A \o c_{(2)}) \o \g^{-1}(d))_{[0]}) \o c_{(1)}\\
&~~~~&\cdot(r_A(id_A \o \varepsilon_C)\nu_{A \o C}((1_A \o c_{(2)}) \o \g^{-1}(d)))_{[1]} \\
&\stackrel{(4.4)}{=}& \b(r_A(id_A \o \varepsilon_C)\xi\nu_{C \o A}(( c_{(2)} \o 1_A) \o \g^{-1}(d))_{[0]}) \o c_{(1)}\\
&~~~~&\cdot(r_A(id_A \o \varepsilon_C)\xi\nu_{C \o A}(( c_{(2)} \o 1_A) \o \g^{-1}(d))_{[1]})\\
&\stackrel{(4.5)}{=}& \sum_{i} \b({r_A(id_A \o \varepsilon_C)\xi({p_i}_{(2)} \o \b^{-1}(b_i))}_{[0]}) \o {p_i}_{(1)}({r_A(id_A \o \varepsilon_C)\xi({p_i}_{(2)} \o \b^{-1}(b_i))}_{[1]})\\
&=&\sum_{i} \b(r_A(id_A \o \varepsilon_C)({b_i}_{[0]} \o \g^{-1}({p_i}_{(2)}){b_i}_{[1]})_{[0]}) \o {p_i}_{(1)}(r_A(id_A \o \varepsilon_C)({b_i}_{[0]} \o \g^{-1}({p_i}_{(2)}){b_i}_{[1]})_{[1]})\\
&=& \sum_{i}\b({b_i}_{[0]}) \o {p_i}_{(1)}\varepsilon_C({p_i}_{(2)})({b_i}_{[1]})\\
&=& \sum_{i}\xi(p_i \o b_i)
= \xi(\nu_{C \o A}((h \o 1_A) \o g)).
\end{eqnarray*}
Use the fact that $\nu_{A \o C}$ is a morphism of right $(C, \g)$-Hom comodules,  we also have
\begin{eqnarray*}
&&\theta(\a^{-1}(c) \o d_{(1)}) \o \g(d_{(2)})\\
&=& r_A(id_A \o \varepsilon_C)\nu_{A \o C}((1_A \o \g^{-1}(c))\o d_{(1)})\o \g(d_{(2)})\\
&=& \sum_{i}r_A(id_A \o \varepsilon_C)(\b^{-1}(a_i) \o {q_i}_{(1)}) \o \g({q_i}_{(2)})\\
&=& \sum_{i} a_i \o q_i = \nu_{A \o C}((1_A \o c) \o d)\\
&\stackrel{(4.4)}{=}& \xi(\nu_{C \o A}((c \o 1_A) \o d)).
\end{eqnarray*}
Hence, we can get condition Eq.(4.1).

 \section{\textbf{Applications}}
\def\theequation{5. \arabic{equation}}
\setcounter{equation} {0}

\subsection {A Maschke type theorem for Doi Hom-Hopf modules}
Since separable functors reflect well the semisimplicity of the objects of a categogy, by Theroem 4.2,  we will prove a Maschke type theorem for   Doi Hom-Hopf modules as an application.\\

\noindent{\bf Corollary 5.1.} Let $(H, A, C)$ be a  Doi Hom-Hopf datum, and
  $M, N \in \widetilde{\mathscr{H}}(\mathscr{M}_k)(H)^{C}_{A}$.
   Suppose that there exists
  a normalized $(A, \b)$-integral
$
  \theta:   (C, \g)  \otimes  (C, \g) \rightarrow  (A, \b)
$. Then a monomorphism ({\it resp.} epimorphism)
  $f : (M, \mu)  \rightarrow  (N, \nu) $
  splits in $\widetilde{\mathscr{H}}(\mathscr{M}_k)(H)^{C}_{A}$, if the monomorphism ({\it resp.} epimorphism) $f  $ splits as an $(A, \b)$-module morphism.

\subsection { Relative Hom-Hopf modules}

Let $(H, \a)$ be a monoidal Hom-Hopf algebra and $(A, \b)$ a  right $(H, \a)$-Hom-comodule algebra. Then the triple $(H, A, H)$ is a  Doi-Hom Hopf datum.
The category $\widetilde{\mathscr{H}}(\mathscr{M}_k)(H)^{H}_{A}$ is called an  $(H,A)$-Hom-Hopf module category and denoted by $\widetilde{\mathscr{H}}(\mathscr{M}_k)^{H}_{A}$.\\

\noindent{\bf Corollary 5.2.}
Let $(H, \a)$ be a monoidal Hom-Hopf algebra and $(A, \b)$ a  right $(H, \a)$-Hom-comodule algebra. Then the following statements are equivalent:

(1) The forgetful functor $F: \widetilde{\mathscr{H}}(\mathscr{M}_k)^{H}_{A}\rightarrow \widetilde{\mathscr{H}}(\mathscr{M}_k)_{A}$ is separable,

(2) There exists a normalized $(A, \b)$-integral $\theta: (H, \a) \otimes (H, \a)  \rightarrow (A, \b)$.

\subsection { Hom-Yetter-Drinfeld modules}

First, we give the definition of Yetter-Drinfeld modules over a monoidal Hom-Hopf algebra,
which is also introduced by Liu and Shen in \cite{LS} similarly.\\

\noindent{\bf Definition 5.3.} Let $(H, \a)$ be a monoidal Hom-Hopf algebra. A right-right $(H, \a)$-Hom-Yetter-
Drinfeld module is an object $(M, \b)$ in $\widetilde{\mathscr{H}}(\mathscr{M}_k)$, such that $(M, \b)$ a right $(H, \a)$-Hom-module  and a  right $(H, \a)$-Hom-comodule  with the following compatibility condition:
\begin{eqnarray}
m_{[0]}\c h_{(1)}\o m_{[1]}\c h_{(2)}=\mu((\mu^{-1}(m)\c h_{(2)})_{[0]})\o h_{(1)}(\mu^{-1}(m)\c h_{(2)})_{[1]}
\end{eqnarray}
for all $h\in H$ and $m\in M$. We denote by $\mathscr{HYD}_H^H$ the category of right-right $(H, \a)$-Hom-Yetter-Drinfeld
modules, morphisms being right $(H, \a)$-linear right $(H, \a)$-colinear maps.
\\

\noindent{\bf Poposition 5.4.}  One has that Eq. (5.1) is equivalent to the following equation:
\begin{eqnarray*}
\rho(m\c h)=m_{[0]}\c \a(h_{(2)(1)})\o S(h_{(1)})(\a^{-1}(m_{[1]})h_{(2)(2)}),
\end{eqnarray*}
for all $h\in H$ and $m\in M$.

{\it Proof.} Similar to \cite{LS}.\\

\noindent{\bf Theorem 5.5.} Let $(H, \a)$ be a monoidal Hom-Hopf algebra with bijective antipode.

(1) $H$ can be made into a right $H^{op} \o H$-Hom-comodule algebra. The coaction $H\rightarrow H\o (H^{op}\o H)$ is given by the formula
\begin{eqnarray*}
h\mapsto  \a(h_{(2)(1)})\o (h_{(2)(2)}\o S(\a^{-1}(h_{(1)}))).
\end{eqnarray*}

(2) $H$ can be made into a right  $H^{op} \o H$-Hom module coalgebra. The action of $H^{op} \o H$ on $H$ is given by the formula
\begin{eqnarray*}
c\triangleleft (h\o k)=\a(k)(\a^{-1}(c)h).
\end{eqnarray*}

(3) The category $\mathscr{HYD}_H^H$  of right-right Hom-Yetter-Drinfeld modules is isomorphic to a
category of Doi Hom-Hopf modules, namely $\widetilde{\mathscr{H}}(\mathscr{M}_k)(H^{op} \o H)^{H}_{H}$.

{\it Proof.} (1) Let us first prove that $H$ is a right $H^{op} \o H$-Hom comodule. For all $h\in H$,
\begin{eqnarray*}
&&(\a^{-1}\o \D_{H^{op} \o H})\rho_H(h)\\
&=&h_{(2)(1)}\o \D_{H^{op} \o H}( h_{(2)(2)}\o S(\a^{-1}(h_{(1)})))\\
&=& h_{(2)(1)}\o  h_{(2)(2)(1)}\o S(\a^{-1}(h_{(1)(2)}))\o h_{(2)(2)(2)} \o  S^{-1}(\a^{-1}(h_{(1)(1)})) \\
&=& \a(h_{(2)(1)(1)})\o h_{(2)(1)(2)}\o S(\a^{-1}(h_{(1)(2)}))\o   \a^{-1}(h_{(2)(2)}) \o  S(\a^{-1}(h_{(1)(1)})) \\
&=& \a^2(h_{(2)(2)(1)(1)})\o \a(h_{(2)(2)(1)(2)})\o S(\a^{-1}(h_{(2)(1)}))\o  h_{(2)(2)(2)} \o S(\a^{-2}(h_{(1)}))\\
&=& \a^2(h_{(2)(1)(2)(1)})\o  \a(h_{(2)(1)(2)(2)})\o S(h_{(2)(1)(1)})\o \a^{-1}(h_{(2)(2)}) \o  S^{-1}(\a^{-2}(h_{(1)})) \\
&=&\rho(\a(h_{(2)(1)}))\o \a^{-1}(h_{(2)(2)})\o S(\a^{-2}(h_{(1)}))\\
&=&(\rho_{H}\o \a^{-1})\rho_H(h)
\end{eqnarray*}
and therefore $H$ is a right $H^{op} \o H$-Hom comodule.

We also have that
\begin{eqnarray*}
\rho(hg)&=&\a(h_{(2)(1)}g_{(2)(1)})\o (h_{(2)(2)}g_{(2)(2)}\o S(\a^{-1}(h_{(1)}g_{(1)}))\o )\\
&=&\a(h_{(2)(1)})\a(g_{(2)(1)})\o (g_{(2)(2)}h_{(2)(2)}\o S(\a^{-1}(g_{(1)}))S(\a^{-1}(h_{(1)})) )\\
&=&(\a(h_{(2)(1)})\o (h_{(2)(2)}\o S(\a^{-1}(h_{(1)}))))(\a(g_{(2)(1)})\o ( g_{(2)(2)}\o S(\a^{-1}(g_{(1)}))))\\
&=& \rho_{H}(h)\rho_{H}(g)
\end{eqnarray*}
(2) We will first prove that $H$ is a $H^{op} \o H$-Hom comodule. For all $h, l, k, m, c\in H$, we have that
\begin{eqnarray*}
[c\triangleleft(h\o k)]\triangleleft(\a(l)\o \a(m))&=&[\a(k)(\a^{-1}(c)h)]\triangleleft(\a(l)\o \a(m))\\
&=&\a^{2}(m)[[k(\a^{-2}(c)\a^{-1}(h))]\a(l)]\\
&=&\a^{2}(m)[[(\a^{-1}(k)\a^{-2}(c))h]\a(l)]\\
&=&\a^{2}(m)[(k\a^{-1}(c))(hl)]\\
&=&[c(mk)](\a(h) \a(l))\\
&=&\a(c)\triangleleft(hl\o mk)=\a(c)\triangleleft[(h\o k)(l\o m)],
\end{eqnarray*}
and this implies that $H$ is a $H^{op} \o H$-Hom comodule.

Using the fact that $(H, \a)$ is an $(H, \a)$-Hom-bimodule algebra, we can obtain that
that $(H, \a)$ is a left $H^{op} \o H$-Hom-module coalgebra.

(3) Let $(M, \mu, \c, \rho_M)$ be such that $(M, \c)$ is a right $(H, \a)$-module and $(M, \rho_M)$ is a right
$(H, \a)$-comodule. Then $M\in \widetilde{\mathscr{H}}(\mathscr{M}_k)(H^{op} \o H)^{H}_{H}$ if and only if
\begin{eqnarray*}
\rho_M(m\c h)&=&m_{[0]}\c \a(h_{(2)(1)})\o m_{[1]}\triangleleft(h_{(2)(2)}\o S(\a^{-1}(h_{(1)})))\\
&=& m_{[0]}\c \a(h_{(2)(1)})\o S(h_{(1)}) (\a^{-1}(m_{[1]})h_{(2)(2)})
\end{eqnarray*}
for all $h\in H$ and $m\in M$. This shows that $\widetilde{\mathscr{H}}(\mathscr{M}_k)(H^{op} \o H)^{H}_{H}$ is isomorphic to $\mathscr{HYD}_H^H$.

From Theorem 4.2 we obtain immediately the following version of
Maschke type theorem for Hom Yetter-Drinfel'd modules:

\noindent{\bf Corollary 5.6.} Let $(H, \a)$ be a monoidal Hom-Hopf algebra with bijective antipode. Then the following
statements are equivalent:

(1) The forgetful functor $F: \mathscr{HYD}_H^H \rightarrow \widetilde{\mathscr{H}}(\mathscr{M}_k)_{H}$ is separable,

(2) There exists a normalized $(H, \a)$-integral $\theta: (H, \a) \otimes (H, \a)  \rightarrow (H, \a)$.

\subsection{ Doi Hom-Hopf Datum $(k, k,H)$ }

Let $(H, \a)$ be a finite-dimensional monoidal Hom-
Hopf algebra.  Recall
from that \cite{CWZ13} that $\varphi\in H^\ast$ is called a right integral on $H^\ast$ if $\varphi \ast h^\ast  =<h^\ast, 1_H>\varphi$ and $\a^\ast\varphi=\varphi$ for all $h^\ast \in H^\ast$, or, equivalently, if $\varphi: H\rightarrow k$  is right
$(H, \a)$-colinear. Now suppose that there exists a normalized $k$-integral $\theta: (H, \a) \otimes (H, \a)  \rightarrow k$. The map $\varphi$ defined by $<\varphi, h>=\theta(1_H\o h)$ is a right integral. Conversely, if $\varphi\in \int^l_{H\ast}$, the $k$-module consisting of classical
integrals on $(H, \a)$, then $\theta(h\o g)=\varphi(gS^{-1}(h))$ is a $k$-integral. This can be proved directly.\\

\noindent{\bf Corollary 5.7.} With the notations as above, then the following statements are equivalent:

(1) The forgetful functor $F: \widetilde{\mathscr{H}}(\mathscr{M}_k)^{H}\rightarrow \widetilde{\mathscr{H}}(\mathscr{M}_k)_{k}$(the category of all vector spaces) is separable.

(2)There exists a  normalized $k$-integral $\theta: (H, \a) \otimes (H, \a) \rightarrow k$ such that the following conditions are satisfied:
\begin{equation*}
\a(h _{(2)}) \theta (\a^{-1}(g) \otimes h _{(1)})=g _{(1)} \theta (g _{(2)} \otimes \a^{-1}(h)  ).
\end{equation*}
\begin{equation*}
\theta(h _{(1)} \otimes h _{(2)})=\varepsilon _{H}(h).
\end{equation*}

\noindent{\bf Corollary 5.8.}
Let $H$ be a finite dimensional cosemisimple monoidal Hom-Hopf algebra. The forgetful functor
$F: \widetilde{\mathscr{H}}(\mathscr{M}_k)^{H}\rightarrow \widetilde{\mathscr{H}}(\mathscr{M}_k)_{k}$ is separable.

\section*{Acknowledgements}

 The work is supported by the NSF of Jiangsu Province (No. BK2012736) and  the Fund of Science and Technology Department of Guizhou Province (No. 2014GZ81365).

\end{document}